\documentclass{article}
\usepackage{graphicx,epstopdf,subfigure}
\usepackage{amsmath}
\usepackage{amsfonts}

\numberwithin{equation}{section}

\begin{document}

\title{Reaction-Diffusion Systems in Epidemiology}
\author{Sebastian Ani\c{t}a$^a$   and Vincenzo Capasso$^{b}$}
\date{$^a$Faculty of Mathematics, ``Alexandru Ioan Cuza'' University of Ia\c{s}i,\\ and ``Octav Mayer'' Institute of Mathematics of the Romanian Academy,\\ Ia\c{s}i 700506,
Romania\\
 {\small email:  sanita@uaic.ro} \\
$^b$ADAMSS \\ (Centre  for Advanced Applied Mathematical and Statistical Sciences)\\
University of Milan\\
Via Saldini 50, 20133 Milano, Italy \\
{\small email:  vincenzo.capasso@unimi.it}}


\maketitle

\centerline{\date{\today}}

\begin{abstract}
A key problem in modelling the evolution dynamics of infectious
diseases is the mathematical representation of the mechanism of
transmission of the contagion which depends upon the way the
specific disease is communicated among different populations or
subpopulations. Compartmental models   describing  a finite number
of subpopulations  can be described mathematically via systems of
ordinary differential equations. The  same  is not possible for
populations which exhibit some continuous structure, such as space
location, age, etc. In particular when dealing with populations
with space structure the relevant quantities are spatial
densities,  whose evolution in time  requires now nonlinear
partial differential  equations, which are known as
reaction-diffusion systems. In this  chapter  we   are presenting
an (historical)  outline of mathematical epidemiology, paying
particular attention to the role of spatial heterogeneity and
dispersal in the population dynamics of infectious diseases. Two
specific examples are discussed,  which  have been the subject of
intensive research by the authors of the present chapter, i.e.
man-environment-man epidemics, and malaria. In addition to the
epidemiological relevance of these epidemics all over the world,
their treatment  requires  a large amount of different
sophisticate  mathematical methods, and has  even posed new non
trivial mathematical  problems, as one can  realize from the list
of references. One of the most relevant problems   faced  by the
present authors, i.e. regional control, has   been emphasized
here: the public health concern consists of eradicating  the
disease in the relevant population, as fast as possible. On the
other hand, very often the entire domain of interest for the
epidemic, is either unknown, or difficult to manage  for an
affordable implementation of suitable environmental sanitation
programmes. This is the reason why regional control has been
proposed;  it might be sufficient to implement such programmes
only in a given subregion conveniently chosen so to lead to an
effective (exponentially fast) eradication of the epidemic in the
whole habitat; it is evident that this practice may have an
enormous importance in real cases with respect to both financial
and practical affordability.
\end{abstract}


\section{Introduction}

Apart from D. Bernoulli (1760) \cite{Bernoulli1760}, those who
established the roots of this field of research (in chronological
order) were:  W. Farr (1840) \cite{Farr1840}, W.H. Hamer (1906)
\cite{Hamer06}, J. Brownlee (1911) \cite{Brownlee11}, R. Ross
(1911) \cite{ross}, E. Martini (1921) \cite{Martini_1921}, A.
J. Lotka (1923) \cite{Lotka23}, W.O. Kermack and A. G. McKendrick
(1927) \cite{Kermack_McKendrick27}, H. E. Soper (1929)
\cite{Soper29}, L. J. Reed and W. H. Frost (1930)
\cite{Frost_1976}, \cite{Abbey52}, M. Puma (1939)
\cite{Puma_1939}, E. B. Wilson and J. Worcester (1945)
\cite{WW_PNAS_1945}, M. S. Bartlett (1949) \cite{Bartlett49}, G.
MacDonald (1950) \cite{Macdonald50}, N.T.J. Bailey (1950)
\cite{Bailey50}, before many others; the  pioneer work by  En'ko
(1989) \cite{Enko1989} suffered  from being written in Russian;
historical accounts of epidemic theory can be found in
\cite{Serfling_1952}, \cite{Dietz_1997}, \cite{Dietz_2009}. After
the late $'70$'s there has been an explosion of interest in
mathematical epidemiology, also thanks to the establishment of a
number  of new journals dedicated to mathematical  biology.

The scheme of this presentation is the following: in Section
\ref{compartmental_models} a general structure of mathematical
models for epidemic systems is presented in the form of
compartmental systems;   in Section \ref{structured_populations}
the concept of field of forces of infection is discussed for
structured populations.


\subsection{Compartmental models} \label{compartmental_models}

\begin{figure}
\setlength{\unitlength}{1cm}
\begin{picture}(20,4)
\put(2,2.5){\vector(0,-1){0.5}}
\put(1.0,2.85){\textit{births
without}} \put(1.0,2.55){\textit{immunity}}
\put(1.5,1){\framebox(1,1){S}}
\put(2.5,1.6){\vector(1,0){2}}
\put(2,1){\vector(0,-1){0.5}}
\put(1.4,0.25){\textit{deaths}}


\put(4.5,1){\framebox(1,1){E}}
\put(5.5,1.6){\vector(1,0){2}}
\put(5,1){\vector(0,-1){0.5}}
\put(4.4,0.25){\textit{deaths}}

\put(7.5,1){\framebox(1,1){I}}
\put(8.5,1.6){\vector(1,0){2}}
\put(8,1){\vector(0,-1){0.5}}
\put(7.4,0.25){\textit{deaths}}

\put(10.5,1){\framebox(1,1){R}}
\put(11,1){\vector(0,-1){0.5}}
\put(10.4,0.25){\textit{deaths}}

\end{picture}
\caption{The transfer diagram for an SEIR compartmental model
including  the susceptible class S, the exposed, but not yet
infective,  class E, the infective class I, and the removed class
R} \label{Fig399}
\end{figure}
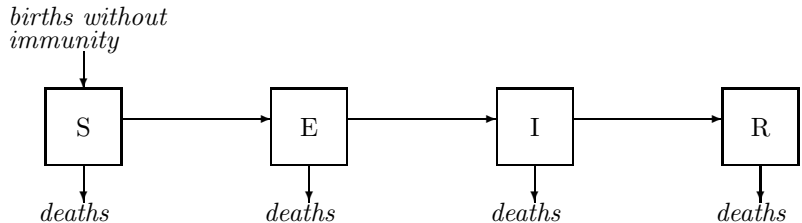

Model reduction  for epidemic systems is obtained via the
so-called compartmental models. In a compartmental model the total
population (relevant to the epidemic process) is divided into a
number (usually small) of discrete categories: susceptibles,
infected but not yet infective (latent), infective, recovered and
immune, without distinguishing different degrees of intensity of
infection; possible structures in the relevant population can be
superimposed when required (see e.g. Figure \ref{Fig399}).

A key problem in modelling the evolution dynamics of infectious
diseases is the mathematical representation of the mechanism of
transmission of the contagion which depends upon the way the
specific disease is communicated among different populations or
subpopulations. This problem has been raised since the very first
models when age and/or space dependence had to be taken into
account.

Suppose at first that the population in each compartment does not
exhibit any structure (space location, age, etc.). Let us ignore,
for the time being,  the intermediate state $E.$ The infection
process ($S$ to $I$) is driven by a force of infection ($f.i.$)
due to the pathogen material produced by the infective population
and available at time $t$

$$(f.i.)(t) = [g(I(.))] (t) $$
which acts upon each individual in the susceptible class. Thus a
typical rate of the infection process is given by the

$$ (\textnormal{incidence rate})(t) = (f.i.)(t) S(t).$$

From this point of view, the so called  ``law of mass action"
simply corresponds to choosing a linear dependence of $g(I)$  upon
$I$

$$(f.i.)(t) = k I(t).$$
The  great  advantage, from a mathematical point of view, is that
the evolution of the epidemic is described (in the space and time
homogeneous cases) by systems of ODE 's which contain at most
bilinear terms.

Referring  to the ``law of mass action'', Wilson  and Worcester \cite{WW_PNAS_1945}
stated the following:

``It would in fact be  remarkable,  in a situation so complex as
that of the passage of an epidemic over a community,  if any
simple law adequately represented the phenomenon in detail ...
even to assume that  the new case rate should be set equal to any
function  ...  might be questioned''.

Indeed  Wilson  and Worcester \cite{WW_PNAS_1945}, and   Severo
\cite{Severo_1969}  had been among the first epidemic modelers
including nonlinear forces of infection of the form

$$(f.i.)(t)= \kappa I(t)^p S(t)^q $$
in their investigations. Here  $I(t)$ denotes the number  of
persons who are infective, and  $S(t)$ denotes the number of
persons who are susceptible to the infection.

Independently, during the analysis  of data regarding the spread
of a cholera epidemic in Southern Italy during 1973, in \cite{18.}
one of the authors (V.C.) suggested the need  to introduce a
nonlinear force of infection  in order  to explain the specific
behavior emerging from the available data.

A  more extended  analysis for a variety of proposed
generalizations of the classical models known as
Kermack-McKendrick  models, appeared  in \cite{12.}, though
nonlinear models  became widely accepted in the literature only a
decade  later,  after  the paper  \cite{Liu_Hethcote_Levin_1987}.

Nowadays models with nonlinear   forces of infection are  analyzed
within the study of various kinds of diseases;  typical
expressions include the so called Holling type functional
responses (see e.g. \cite{12.},
\cite{Hethcote_vandendriessche_91})

$$(f.i.)(t)=g(I(t))\,\,\,;$$
with

\begin{equation} \label{nonlinear_2}
g(I)=\displaystyle{{k\,\, I^p}\over{\alpha +
\beta\,\,\displaystyle{ I^q}}}\,\,\,\,,\qquad p, q\,\, >0
\,\,\,.\end{equation}

Particular cases are

\begin{equation} \label{nonlinear_1}
g(I)=k\,\, I^p\,\,\,\,, \qquad p>0
\end{equation}

\begin{figure}
\centering\includegraphics[width=12.0cm]{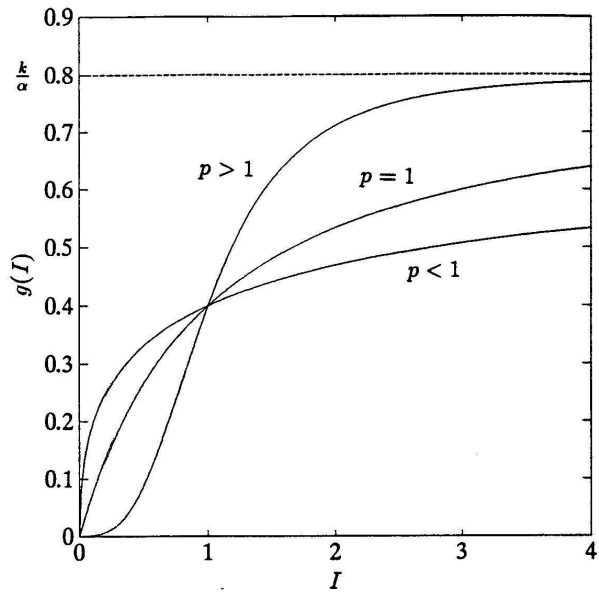}
 \caption{Nonlinear forces of infection \cite{12.}} \label{Fig3.4}
\end{figure}

For the  case  $p=q$  we have the behaviors described in  Figure
\ref{Fig3.4}.

Additional shapes of  $g(I),$  as proposed  in \cite{12.} which
may decrease for large values of  $I,$  may be interpreted as
``awareness'' effects in the contact rates. Significant
contributions to this concept and related epidemiological issues
in recent literature  can be found in \cite{DOnofrio}.

Further extensions  include a nonlinear dependence upon both $I$
and $S$, as discussed in  modelling  AIDS epidemics (see e.g.
\cite{Castillo_Chavez_1_1989}, \cite{Castillo_Chavez_2_1989}, and
references therein), where the social structure of the host
population is analyzed too.


\subsection{Structured populations} \label{structured_populations}

Compartmental models with a finite number of compartments can be
described mathematically via systems of ordinary differential
equations (ODE's). The  same  is not possible for  populations
which exhibit some continuous structure (identified here by a
parameter $z$),   such as  space location, age, etc.

When dealing with populations with space structure the relevant
quantities are spatial densities, such as $s(z;t)$ and $i(z;t)$,
the spatial densities of susceptibles and of infectives
respectively, at a point $z$ of the habitat $\Omega$, and at time
$t \geq 0,$ such that the  corresponding total populations are
given by

$$S(t)=\int_\Omega {s(z;t)\,\, dz}, \qquad  I(t)=\int_\Omega {i(z;t)\,\, dz }$$

The role of spatial heterogeneity and dispersal in population
dynamics has been the subject of much research. We shall refer
here to the fundamental literature on the subject by quoting
\cite{aris_75}, \cite{aronson_weinberger}, \cite{fife},
\cite{Kendall_57}, \cite{levin_78}, \cite{murray_89},
\cite{okubo}, \cite{skellam}, \cite{smoller}.\par In the theory of
propagation of infectious diseases the motivation of the study of
the effects of spatial diffusion is mainly due to the large scale
impact of an epidemic phenomenon; this is discussed in detail in
\cite{cliff_ord}.\par We will choose here as a good starting point
the pioneer work by D.G. Kendall who modified the basic
Kermack-McKendrick SIR model to include the effects of spatial
heterogeneity in an epidemic system \cite{Kendall_57}.\par Kendall
's work has motivated a lot of research in the theory of epidemics
with spatial diffusion. In particular two classes of problems
arise according to the size of the spatial domain or habitat.\par
If the habitat is unbounded, travelling waves are of interest
\cite{Kendall_57}, \cite{Hoppensteadt_75}, \cite{Aronson_1977},
\cite{diekmann_78}, \cite{thieme_77}. A nice introduction to the
subject can be found in \cite{murray_89}.\par If, on the other
hand, the habitat is a bounded spatial domain then problems of
existence of nontrivial endemic states (possibly with spatial
patterns) are of interest \cite{VK78}. \par Usually reaction
diffusion systems (see e.g. \cite{fife}) are seen as an extension
of compartmental systems in which each compartment, representing a
different species, is allowed to invade a spatial domain $\,
\Omega \subset \mathbb R^m \,$ with a space dependent density.
Densities interact among themselves according to the same
mathematical laws which were used for the space independent case,
but are subject individually to a spatial diffusion mechanism
usually committed to the Laplace operator, which simulates random
walk or Brownian motion of the interacting species
\cite{aronson_85}, \cite{fife}.\par

Typically then a system of $\,n\,$ interacting species, each of
them having a spatial density
$$\left\{u_i(x;t),\quad x\in \Omega \right\}, \quad i=1,\ldots,n \quad
\hbox{\rm at time} \quad t \geq0,$$ \vskip 8pt

\noindent is
described by the following system of semilinear parabolic
equations: \vskip 8pt

\begin{equation} \label{reaction_diffusion}
\displaystyle{\frac{\partial u}{\partial t}=D \Delta
u(x;t)+f(u(x;t))} \end{equation}
 \vskip 8pt

\noindent in $\Omega \times \mathbb R_+ ,$ subject to suitable boundary
conditions.\par \noindent Here $\,u=(u_1,\ldots,u_n)^T\, ;\,\,\,
D=diag(d_1,\ldots,d_n),\,\,\, \hbox{\rm and} \,\,\, f(z)\,,\,\,\,
z\in \mathbb R^n\,$ is the interaction law among the species via
their densities.\par A derivation and discussion of Equation
(\ref{reaction_diffusion}) can be found in \cite{fife} (see also
\cite{okubo}).\par In Equation (\ref{reaction_diffusion})
$\,f(z)\,,\,\,\, z\in \mathbb R^n, $ usually is the same
interaction function of the classical compartmental approach. Thus
for an SIR model (including Susceptible, Infective, and Removed
individuals) with vital dynamics, Equation
(\ref{reaction_diffusion}) is  written as follows

\begin{equation}  \label{5.2}
\left\{ \begin{array}{ll} \displaystyle{\frac{\partial s}
{\partial t}(x;t)=d_1 \Delta s(x;t)- k\,i(x;t)\,s(x;t)+\mu-\mu\,
s(x;t)} \\ \smallskip \\\displaystyle{\frac{\partial i}{\partial
t}(x;t) = d_2 \Delta
i(x;t)+k\,i(x;t)\,s(x;t)-(\mu+\gamma)\,i(x;t)} \\ \smallskip \\
\displaystyle{\frac{\partial r}{\partial t}(x;t)= d_3 \Delta
r(x;t)+\gamma\, i(x;t)- \mu \, r(x;t)}
 \end{array} \right.
\end{equation}

\noindent where now $\,
s(x;t)\,,\,\,i(x;t)\,,\,\,r(x;t)\,$ are the corresponding spatial
densities of $S, I, R$:
$$S(t)=\int_{\Omega} s(x;t)\, dx; \quad I(t)=\int_{\Omega} i(x;t)\, dx;\quad
R(t)=\int_{\Omega} r(x;t)\, dx. $$ \vskip 8pt
Thus the infection
process is represented by a \hskip 4pt ``local'' interaction of two
densities $s(x;t)$ and $i(x;t),$ at point $x \in \Omega$ and time
$t \in \mathbb R_+,$ via the following ``law of mass action''

\begin{equation} \label{5.3} (f.i.)(z;t)= k \,i(x;t)\,s(x;t) \end{equation}

In Kendall 's model \cite{Kendall_57} (\ref{5.3}) is
substituted by an ``integral'' interaction between the infectives
and the susceptibles; the force of infection acting at point $x\in
\Omega$ and time $t\in \mathbb R_+$ is given by
\begin{equation} \label{5.4}
g(i(\cdot;t))(x)= \int_{\Omega} k(x,x')\, i(x';t)\, dx',
\end{equation}
 where $k(x,x')$ describes the influence, by any reason, of the
infectives located at any point $x' \in \Omega$ on the
susceptibles located at point $x \in \Omega$.

As a consequence the infection process is described now by
\begin{equation} \label{5.5}
(f.i.)(z;t)= g(i(\cdot;t))(x)\, s(x;t).
\end{equation}

 Clearly we reobtain (\ref{5.3}) in the limiting case $\,
k(x,x')=\delta (x-x')\,$ (the Dirac function).

As one can see, (\ref{5.5}) better fits the philosophy introduced
in Section 1 to allow $g(I)$, the force of infection due to the
population of infectives, and acting on the susceptible
population, to have a general form which case by case takes into
account the possible mechanisms of transmission of the
disease.\par For Problem (\ref{5.2}) we refer to the literature
\cite{Aronson_1977}, \cite{VK78}, \cite{diekmann_78},
\cite{thieme_77}. For the   case in (\ref{5.4})-(\ref{5.5}) the
emergence of travelling waves has been shown
 in \cite{Kendall_57} and \cite{Aronson_1977}. The analysis of the diffusion
approximation of Kendall's model can be found in
\cite{Hoppensteadt_75}.

When dealing with populations with an age structure,  we may
interpret the parameter $z$ as the age-parameter so that the first
 model above  is a model with intracohort interactions while the second one
 is a model with intercohort interactions  (see  e.g. \cite{Busenberg_etal_1988}, and references therein).

A  large literature on the subject can be found in \cite{VK09}.

\section{Specific examples}

\subsection{Spatially  structured  man-environment-man  epidemics}
\label{spatial_structure}

A  widely  accepted model for  the spatial spread of epidemics in
an habitat $\Omega $, via the environmental pollution  produced by
the infective population, e.g.  via  the excretion of pathogens in
the environment, is the following  one, as proposed in
\cite{VK78}, \cite{VK84} (see also \cite{VK09}, and references
therein). The  model below  is a more realistic generalization of
a previous model proposed by one of the authors (V. C.) and his
co-workers \cite{VK81}, \cite{7}, \cite{18.} to describe
fecal-orally transmitted diseases (cholera, typhoid fever,
infectious hepatitis, etc.) which are typical for the European
Mediterranean regions; it can anyhow be applied to other
infections, and other regions,  which are propagated by similar
mechanisms (see  e.g. \cite{Cod}); schistosomiasis in Africa is a
typical additional example \cite{nasell_hirsch}.

\begin{equation}  \label{(1.1)}
\left\{ \begin{array}{ll} \displaystyle {{\partial
u_1}\over{\partial t}}(x,t)=d_1\Delta u_1(x,t) \displaystyle
-a_{11}u_1(x,t)+\int_{\Omega }k(x,x') u_2(x',t)dx'
\vspace{2mm}\\
\displaystyle {{\partial u_2}\over{\partial t}}(x,t)=
-a_{22}u_2(x,t)+g(u_1(x,t))
          \end{array}
 \right.  \end{equation}
 in  $\Omega \subset {\bf R}^N\, \ (N\geq 1)$, a nonempty bounded
domain with a smooth boundary $\partial \Omega $; for  $t \in
(0,+\infty )$, where $a_{11}\geq 0$, $a_{22}\geq 0,~d_1>0$ are
constants.

\begin{itemize}
\item[$\bullet$] $u_1(x,t)$ denotes the concentration of the
pollutant (pathogen material)  at a  spatial location $x \in
\overline{\Omega },$ and a time $t\geq 0;$

\item[$\bullet$] $u_2(x,t)$ denotes the spatial distribution of
the infective population.

\item[$\bullet$]The terms  $ -a_{11}u_1(x,t)$ and
$-a_{22}u_2(x,t)$ model natural decays.

\item[$\bullet$] The total susceptible population is assumed to be
sufficiently large with respect to the infective population, so
that it can be taken as constant.

\end{itemize}

For this kind of epidemics the infectious agent is multiplied by
the infective human population and then sent to the sea through
the sewage; because of the peculiar eating habits of the
population of these regions the agent may return via some
diffusion-transport mechanism to any point of the habitat
$\Omega,$ where the infection process is further activated; thus
the integral term

$$\int_{\Omega }k(x,x') u_2(x',t)dx' $$
expresses the fact that the pollution produced at any point $x'
\in \Omega$ of the habitat is made available at any other point $x
\in \Omega$; when dealing with human pollution, this may  be due
to either malfunctioning of the sewage system, or  improper
dispersal of sewage in the habitat. Linearity  of the above
integral operator  is   just a simplifying option.

The Laplace operator takes into account a simplified random
dispersal of the infectious agent in the habitat $\Omega,$ due to
uncontrolled additional causes of dispersion (with a constant
diffusion coefficient to avoid purely technical complications); we
assume that the infective population does not diffuse  (the case
with diffusion would   be here a technical simplification). As
such, System (\ref{(1.1)}) can be adopted as a good model for the
spatial propagation of an infection in agriculture and forests,
too.

 Finally, the  local ``incidence rate'' at point $x  \in \Omega
,$ and time $t\geq 0$, is  given by

$$(i.r.)(x,t)=g(u_1(x,t)),$$ depending
upon the local concentration of the  pollutant.

 The parameters  $a_{11} $ and $a_{22} $ are intrinsic decay parameters of the two populations.

\subsection{Seasonality.}

If we wish to model a  large class of  fecal-oral transmitted
infectious diseases, such as typhoid  fever, infectious hepatitis,
cholera, etc., we may include the possible seasonal variability of
the environmental conditions, and their impact on the habits of
the susceptible population, so that the relevant parameters are
assumed periodic in time, all with the same period $T\in
(0,+\infty )$.

 As a purely technical simplification, we may assume that only the
incidence rate is periodic, and in particular that it can be
expressed as

$$(i.r.)(x,t)=h(t,u_1(x,t))=p(t)g(u_1(x,t)),$$
where $h$, the functional dependence of the incidence rate upon
the concentration of the pollutant, can be chosen as in the time
homogeneous case, with possible behaviors  as shown in Figure
\ref{Fig3.4}.

The explicit time dependence of the incidence rate is given via
the function $p(\cdot),$  which is assumed to be a strictly
positive, continuous and $T-$periodic function of time; i.e. for
any  $t \in \mathbb R,$

$$p(t)=p(t+T).$$

\bf Remark 1. \rm The  results  can be easily
 extended to the case in which  also $a_{11}$, $a_{22}$ and $k$ are $T-$periodic
functions.

In \cite{VK84} the above model was  studied, and sufficient
conditions were given for either the asymptotic extinction of an
epidemic outbreak, or the existence and stability of an  endemic
state; while in \cite{VK83} the periodic case  was additionally
studied, and sufficient conditions were given for either the
asymptotic extinction of an epidemic outbreak, or the existence
and stability of a periodic endemic state with the same period of
the parameters.

\subsection{Saddle point behaviour.}

The choice of $g$ has a strong influence on the dynamical behavior
of system (\ref{(1.1)}). The case in which $g$ is a monotone
increasing function with constant concavity has been analyzed in
an extensive way (see \cite{VK09}, \cite{VK_Kunisch},
\cite{VK82}); concavity leads to the existence (above a parameter
threshold) of exactly one nontrivial endemic state and to its
global asymptotic stability. In order  to   better  clarify the
situation, consider  first  the spatially homogeneous case (ODE
system) associated with system (\ref{(1.1)}); namely

\begin{equation}\left\{ \begin{array}{ll} \displaystyle {{d
z_1}\over{d t}}(t)=\displaystyle -a_{11}z_1(t)+ a_{12} z_2(t)
\vspace{2mm}\\
\displaystyle {{d u_2}\over{d t}}(t)= -a_{22}z_2(t)+g(z_1(t))
          \end{array}
 \right. \label{ODE} \end{equation}

 In \cite{VK82} and
\cite{nasell_hirsch} the bistable case (in which system
(\ref{ODE}) may admit  two nontrivial steady states, one of which
is a saddle point in the phase plane) was obtained by assuming
that the force of infection, as a function of the concentration of
the pollutant, is sigma shaped. In \cite{nasell_hirsch} this shape
had been  obtained as a consequence of the sexual reproductive
behavior of the schistosomes. In \cite{VK82} (see also
\cite{VK81}) the case of fecal-oral transmitted diseases was
considered; an interpretation of the sigma shape of the force of
infection was proposed to model the response of the immune system
to environmental pollution: the probability of infection is
negligible at low concentrations of the pollutant, but increases
with larger concentrations; it then becomes concave and saturates
to some finite level as the concentration of pollutant increases
without limit.

Let  us now refer  to the following simplified form of   System
(\ref{(1.1)}),  where   as  kernel we have taken  $k(x, x')=
a_{12}\delta(x-x'),$

\begin{equation}\left\{ \begin{array}{ll} \displaystyle {{\partial
u_1}\over{\partial t}}(x,t)=d_1\Delta u_1(x,t) \displaystyle
-a_{11}u_1(x,t)+ a_{12} u_2(x,t)
\vspace{2mm}\\
\displaystyle {{\partial u_2}\over{\partial t}}(x,t)=
-a_{22}u_2(x,t)+g(u_1(x,t))
          \end{array}
 \right. \label{(1.1bis)} \end{equation}

The concavity of $g$ induces concavity of its evolution operator,
which, together with the monotonicity induced by the quasi
monotonicity of the reaction terms in (\ref{(1.1bis)}), again
imposes uniqueness of the possible nontrivial endemic state. On
the other hand, in the case where $g$ is sigma shaped,
monotonicity of the solution operator is preserved, but as we have
already observed in the ODE case, uniqueness of nontrivial steady
states is no longer guaranteed. Furthermore, the saddle point
structure of the phase space cannot be easily transferred from the
ODE to the PDE case, as discussed in \cite{VK82}, \cite{VK97}. In
\cite{VK82}, homogeneous Neumann boundary conditions were
analyzed;  in this case nontrivial spatially homogeneous steady
states are still possible. But  when we deal with homogeneous
Dirichlet boundary conditions or general third-type boundary
conditions, nontrivial spatially homogeneous steady states are no
longer allowed. In \cite{VK97} this problem was faced  in more
detail; the steady-state analysis was  carried out and  the
bifurcation pattern of nontrivial solutions to system
(\ref{(1.1bis)}) was determined when subject to homogeneous
Dirichlet boundary conditions. When the diffusivity of the
pollutant is small,  the existence of a narrow bell-shaped steady
state was shown, representing  very likely  a saddle point for the
dynamics of (\ref{(1.1bis)}). Numerical  experiments confirm the
bistable situation:  ``small'' outbreaks stay localized under this
bell-shaped steady state, while ``large'' epidemics tend  to
invade the whole habitat.

\subsection{Boundary feedback}
An interesting problem concerns the case  of boundary feedback of
the pollutant, which has been proposed  in \cite{VK_Kunisch},  and
further  analyzed  in \cite{VK_Thieme};  an optimal control
problem has been later analyzed  in \cite{VK_Barbu_Arnautu}.

In this  case  the reservoir of the pollutant generated by the
human population is spatially separated from the habitat by a
boundary through which the positive feedback occurs. A model of
this kind has been proposed as an extension of  the ODE  model for
fecal-oral transmitted infections in Mediterranean coastal regions
presented  in \cite{18.}.

For this kind of epidemics the infectious agent is multiplied by
the infective human population and then sent to the sea through
the sewage system; because of the peculiar eating habits of the
population of these regions, the agent may return via some
diffusion-transport mechanism to any point of the habitat, where
the infection process is restarted.

The mathematical model is based on the following system of
evolution equations:

\begin{equation} \label{boundary_feedback}
\left\{
\begin{array}{ll}
\displaystyle{\frac{\partial u_1}{\partial t} }(x;t) =\Delta
u_1(x;t)-a_{11}
u_1(x;t)  \nonumber \\
~~~ \\
\displaystyle{\frac{\partial u_2}{\partial t}}(x;t)
=-a_{22}u_2(x;t)+ g(u_1(x;t))
\end{array}
\right.
\end{equation}

\noindent
in $\Omega \times (0,+\infty)$, subject to the following
boundary condition
$${{\partial u_1}\over {\partial \nu}} (x;t) +\alpha u_1(x;t)=\int_{\Omega} k
(x,x') u_2(x';t)\, dx' $$
 on $\partial \Omega \times
(0,+\infty)$, and also subject to suitable initial conditions.

Here $\Delta$ is the usual Laplace operator modelling the random
dispersal of the infectious agent in the habitat; the human
infective population is supposed not to diffuse. As  usual
$a_{11}$ and $a_{22}$ are positive constants.  In the boundary
condition  the left hand side is the general boundary operator
$B:=\displaystyle{{\partial\over{\partial \nu}}+\alpha(\cdot)}$
associated with the Laplace operator; on the right hand side the
integral operator
$$H\left[u_2(\cdot,t)\right](x):= \int _{\Omega} k(x,x') u_2(x';t)\, dx' $$
describes boundary feedback mechanisms, according to which the
infectious agent produced by the human infective population at
time $t>0$, at any point $x' \in \Omega$, is available, via the
transfer kernel $k(x,x')$, at a point $x \in
\partial \Omega$.

 Clearly the boundary $\partial
\Omega$ of the habitat $\Omega$ can  be divided into two disjoint
parts: the sea shore $\Gamma_1$ through which the feedback
mechanism may occur, and $\Gamma_2$ the boundary on the land, at
which we may assume complete isolation.

The parameter $\alpha (x)$ denotes the rate at which the
infectious agent is wasted away from the habitat into the sea
along the sea shore.
 Thus one may well  assume that
$$\alpha(x),\quad k(x,\cdot) =0, \qquad \hbox{\rm for} \,\,\,\, x\in \Gamma_2\,\,\,.$$

A relevant  assumption, of great importance in the control
problems   that  we have been facing later,  is that   the habitat
$\Omega$  is "epidemiologically" connected to its boundary by
requesting that
$$\hbox{\rm for any} \,\,\,\,\, x'\in\Omega \,\,\,\,\,\hbox{\rm there exists some}\,\,\,\,\,
x\in \Gamma_1 \,\,\,\,\, \hbox{\rm such that} \quad k(x,x')>0. $$

This means that from any point of the habitat infective
individuals contribute to polluting at least some point on the
boundary (the sea shore).

In the above model  delays had been  neglected and the feedback
process had been  considered to be linear; various extensions have
been considered in subsequent  literature.

\subsection{Malaria}  \label{malaria}

Malaria is found throughout the tropical and subtropical regions
of the world and causes more than 300 million acute illnesses and
at least one million deaths annually \cite{WHO_UNICEF}, especially
in Africa.  Human malaria is caused by one or a combination of
four species of plasmodia: Plasmodium falciparum, P. vivax, P.
malariae, and P. ovale. The parasites are transmitted through the
bite of infected female mosquitoes of the genus Anopheles.
Mosquitoes can become infected by feeding on the blood of infected
people, and the parasites then undergo another phase of
reproduction in the infected mosquitoes.

The earliest attempt to provide a quantitative understanding of
the dynamics of malaria transmission was that of Ross \cite{ross}.
Ross' models consisted of a few differential equations to describe
changes in densities of susceptible and infected people, and
susceptible and infected mosquitoes. Macdonald \cite{mcdonald}
extended Ross' basic model, analyzed several factors contributing
to malaria transmission, and concluded that ``the least influence
is the size of the mosquito population, upon which the traditional
attack has always been made''.

The work of Macdonald had a very beneficial impact on the
collection, analysis, and interpretation of epidemic data on
malaria infection \cite{molin} and guided the enormous global
malaria-eradication campaign of his era.

The classical Ross-Macdonald model is highly simplified.
Subsequent contributions have been made to extend the
Ross-Macdonald malaria models considering a variety of
epidemiological features of malaria \cite{AronMay}, \cite{Dietz_1}
(see also  \cite{Bac}, \cite{Chitnis_2006}, \cite{Chitnis_2008},
\cite{Dietz}, \cite{Killeen}, \cite{ruan}, \cite{Smith}).

Here we will consider an oversimplified model, concentrating  on a
problem of  possible eradication of a spatially structured malaria
epidemic, by acting on the segregation of the human and the
mosquito populations via e.g. treated bednets (see e.g.
\cite{Sochanta}).

Following recent papers by Ruan et al \cite{ruan}, and Chamchod
and Britton \cite{britton}, we divide the human population into
two classes, susceptible and infectious, whereas the mosquito
population is divided into three classes, susceptible, infectious,
and removed because of death. Suppose that the infection in the
human does not result death or isolation. For the transmission of
the pathogen, it is assumed that a susceptible human can receive
the infection only by contacting with infective mosquitos, and a
susceptible mosquito can receive the infection only from the
infectious human.

For simplicity, assume the total populations of both humans and
mosquitoes are constants and denoted by $H$ and $M$, respectively.
Let $X(t)$ and $Y(t)$ denote the numbers of infected humans and
mosquitoes at time $t$, respectively. Let $a$ be the rate of
biting on humans by a single mosquito (number of bites per unit
time). Then the number of bites on humans per unit time per human
is $\displaystyle{\frac{a}{H}}$. If $b$ is the proportion of
infected bites on humans that produce an infection, the
interaction between the infected mosquitoes $Y(t)$ and the
uninfected humans $H - X(t)$ will produce new infected humans at
a rate of of $\displaystyle{\frac{a}{H}} b ( H - X(t)) Y(t). $
Malaria on humans is an SIS system  for which  human infectives
after recovery go back in the susceptible state; we will denote by
$r$ the per capita rate of recovery in humans so that $1/r$ is the
duration of the disease in humans \cite{VK09}.

Therefore, the equation for the rate of change in the number of
infected humans is
$$\frac{d X}{d t} = -r X(t) + \frac{a}{H} b ( H - X(t)) Y(t) .$$

Similarly, if $\mu$ is the per capita rate of mortality in mosquitoes, so that $%
1/\mu$ is the life expectancy of mosquitoes,  and $c$ is the
transmission efficiency from humans to mosquito, then we have the
equation for the rate of change in the number of infected
mosquitoes:

$$
\frac{d Y}{d t} =- \mu Y(t)+ \frac{a}{H} cX(t) (M - Y(t)) .
$$

Here we wish to consider a spatially structured system, so that we
will refer to spatial densities of infected humans and infected
mosquitoes. Specifically, let $u_{1}(x,t)$ denote the spatial
density of the population of infected mosquitoes at a spatial
location $x\in \overline{\Omega }$ and a time $t\geq 0;$ while
$u_{2}(x,t)$ denotes the spatial distribution of the human
infective population. In accordance with the above, the spatial
density $C(x)$ of the total human population will be assumed
constant in time, so that $C(x) - u_{2}(x,t)$ will provide the
spatial distribution of susceptible humans, at a spatial location
$x\in \overline{\Omega }$ and a time $t\geq 0.$

Hence the ``local incidence'' for humans, at point $x\in \overline{\Omega }$%
, and time $t\geq 0$, is taken of the form
\[
(i.r.)_H(x,t)=g(u_{1}(x,t),u_{2}(x,t))=(C(x)-u_{2}(x,t))h(u_{1}(x,t)),
\]
depending upon the local densities of both populations via a
suitable functional response $h$.

As pointed out in \cite{molin}, seasonality of the aggressivity to
humans by the mosquito population might also be considered in the
functional response $h;$ this has been a topic of the authors
research programme in \cite{VK83}, and \cite{AKB},
\cite{AKBelleni}.

As proposed in \cite{Bac}, we will include spatial diffusion of
the infective mosquito population (with constant diffusion
coefficient to avoid purely technical complications), but we
assume that the human population does not diffuse.

On the other hand, inspired by previous models for spatially
structured man-environment epidemics (see e.g. \cite{VK84}), we
will modify the infection term of mosquitoes
$\displaystyle{\frac{a}{H}} cX(t) (M - Y(t))$ as follows.

Given a habitat $\Omega \subset \mathbf{R}^{N}\,\ (N\geq 2)$
(which is a nonempty bounded domain with a smooth boundary
$\partial \Omega $), we assume that the total susceptible mosquito
population is so large that it can be considered  time and space
independent; further we consider the fact
that the infected mosquitoes at a location $x\in {\Omega }$ and a time $%
t\geq 0$ is due to contagious bites to human infectives at any point $%
x^{\prime}\in \Omega$ of the habitat, within a spatial
neighborhood  of $x$ represented by a suitable probability kernel
$k(x^{\prime},x),$  depending on the specific structure of the
local ecosystem; as a  trivial simplification one may assume
$k(\cdot,x)$ as a Gaussian density centered at $x;$ hence the
``local incidence'' for mosquitoes, at point $x\in {\ \Omega } $,
and time $t\geq 0$, is taken as

\[
(i.r.)_M(x,t)=\int_{\Omega }k(x,x^{\prime})
u_2(x^{\prime},t)dx^{\prime}
\]

Finally, from now on, we will denote by $-a_{11}u_1(x,t)$ the
natural decay of the infected mosquito population, while
$-a_{22}u_2(x,t)$ will denote the recovery rate (back to
susceptibility) of the human infective population.

All the above leads to the following over simplified model for the
spatial spread of malaria epidemics, in which we have ignored the
possible acquired immunity of humans after exposure to the
contagion, the possible differentiation of the mosquito
population, etc. (see e.g. \cite{Killeen}).

$$
\left\{
\begin{array}{ll}
\displaystyle
{{\partial u_{1}}\over {\partial t}}(x,t)=d_{1}\Delta u_{1}(x,t)\displaystyle %
-a_{11}u_{1}(x,t)+\int_{\Omega }k(x,x^{\prime })u_{2}(x^{\prime
},t)dx^{\prime }\vspace{2mm} &  \\
{{\partial u_{2}}\over {\partial t}}(x,t)=-a_{22}u_{2}(x,t)+g(u_{1}(x,t), u_2(x,t)) &
\end{array}
\right.
$$
in $\Omega \times (0,+\infty )$, where $a_{11}>0$,
$a_{22}>0,~d_{1}>0$ are constant.

\section{Regional  Control: \emph{Think Globally, Act Locally}}

Let us now go back to  System (\ref{(1.1)}) in  $\Omega \subset
{\bf R}^N\, \ (N\geq 1)$, a nonempty bounded domain with a smooth
boundary $\partial \Omega $; for  $t \in (0,+\infty )$, where
$a_{11}\geq 0$, $a_{22}\geq 0,~d_1>0$ are constants.

The public health concern consists of  providing methods for the
eradication of the disease in  the relevant population, as fast as
possible. On the other hand, very often the entire domain
$\Omega,$ of interest for the epidemic, is either unknown, or
difficult to manage  for an affordable  implementation of suitable
environmental sanitation programmes. Think of malaria,
schistosomiasis, and alike,  in Africa, Asia, etc.

This has led the second  author,  in a  discussion with Jacques
Louis Lions in 1989, to suggest that it might be sufficient to
implement such programmes only in a given subregion $\omega
\subset \Omega,$ conveniently chosen so to lead to an effective
(exponentially fast) eradication of the epidemic in the whole
habitat $\Omega.$ Though, a satisfactory mathematical treatment of
this issue has been obtained only few years later in \cite{AK02}.
This practice may have an enormous importance in real cases with
respect to both financial and practical affordability. Further,
since  we propose to act on the elimination of  the pollution
only,  this practice means an additional nontrivial social benefit
on the human population, since it  would not be limited in his
social and alimentation habits.

\begin{figure}
\centering\includegraphics[width=12.0cm]{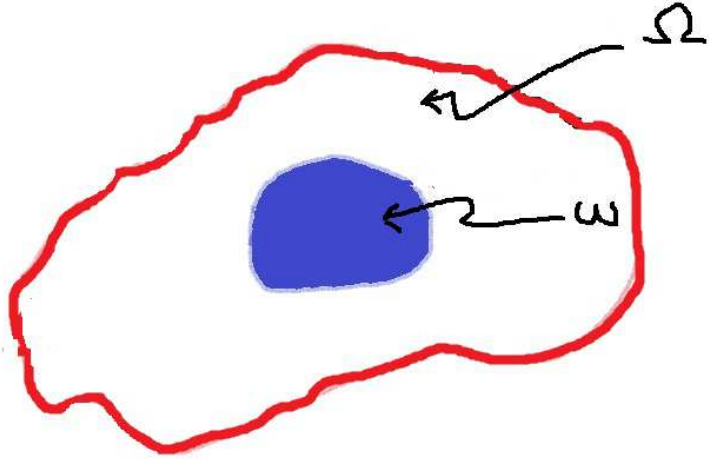}
 \caption{Think Globally, Act Locally.}
\end{figure}\label{Fig2}

In this section   a review is presented of some results obtained
by the authors, during 2002-2012,  concerning stabilization (for
both the time homogeneous case  and  the periodic case).
Conditions have been provided for the exponential decay of the
epidemic in the whole habitat $\Omega,$ based on the elimination
of the pollutant in a subregion $\omega \subset \Omega $. The case
of homogeneous third type boundary conditions  has been
considered, including the homogeneous Neumann boundary conditions
(to mean complete isolation of the habitat):

$$ \frac{\partial u_1}{\partial \nu }(x,t)+\alpha u_1(x,t)=0 \qquad \mbox{\rm on } \partial {\Omega}\times
(0,+\infty),$$

where $\alpha \geq 0$ is a constant.

For the time homogeneous case  the following assumptions  have
been taken:
\begin{itemize}
\item[{\bf (H1)}]   $\, \, g: \mathbb R\rightarrow [0,+\infty )$
is a function satisfying

       \begin{itemize}
      \item[$a)$] $ g(x)=0$, \quad for $x\in (-\infty ,0]$,

       \item[$b)$] $ g$ is Lipschitz continuous and increasing,

      \item[$c)$] $ g(x) \leq a_{21}x, \ \mbox{\rm for any } x\in [0,+\infty )$,
       where $\,\, a_{21}>0 ;$
       \end{itemize}
\medskip
\item[{\bf (H2)}] $k\in L^{\infty }(\Omega \times \Omega ),
 \,\, k(x,x')\geq 0 \,$ a.e. in $\Omega \times \Omega$,
$$\int_{\Omega }k(x,x')dx >0 \qquad \mbox{\rm a.e. } x'\in \Omega ; $$

\item[{\bf (H3)}] $u_1^0, \ u_2^0 \in L^{\infty }(\Omega), \quad
u_1^0(x),
 \ u_2^0(x)\geq0$ \quad a.e. in $\Omega$.
\end{itemize}

Let $\omega\subset \subset \Omega$ be a nonempty subdomain with a
smooth boundary and $\Omega \setminus \overline{\omega }$ a
domain. Denote by $\chi_{\omega}$ the characteristic function of $
\omega$ (we use the convention
$$\chi_{\omega}(x)h(x)=0, \qquad x\in {\bf R}^N \setminus \overline{\omega },$$
even if function $h$ is not defined on the whole set ${\bf R}^N
\setminus \overline{\omega }). $

Our goal is to study the  controlled system
$$\left\{ \begin{array}{ll} \displaystyle {{\partial u_1}\over
{\partial t}}(x,t)= \displaystyle d_1\Delta
u_1(x,t)-a_{11}u_1(x,t)
+\int_{\Omega }&k(x,x') u_2(x',t)dx' \\
\qquad  \qquad \quad +\chi_{\omega}(x)\, v(x,t),\quad &(x,t)\in \Omega \times (0,+\infty ) \\
\displaystyle {{\partial u_1}\over {\partial \nu }}(x,t)+\alpha
u_1(x,t)=0,
\quad &(x,t)\in \partial \Omega \times (0,+\infty ) \\
\displaystyle {{\partial u_2}\over {\partial
t}}(x,t)=-a_{22}u_2(x,t) +g(u_1(x,t)), \ \
&(x,t)\in \Omega \times (0,+\infty ) \\
u_1(x,0)=u_1^0(x), \quad u_2(x,0)=u_2^0(x), &x\in \Omega ,
          \end{array}
  \right. \label{(1.2)}$$
subject to a control   $v\in L^{\infty }_{loc}(\overline{\omega
}\times [0,+\infty ))$ (which implies that $supp (v(t))\subset
{\overline \omega}$ for $t\geq 0$).

We have to mention that existence, uniqueness and nonnegativity of
a solution to the  above system can be proved as in
\cite{babak_2007}. The nonnegativity of $u_1$ and $u_2$ is a
natural requirement due to the biological significance of $u_1$
and $u_2$.

\bf Definition 3.1. \rm
We say that our system  is \emph{zero-stabilizable} if
 for any $u_1^0$ and $u_2^0$ satisfying (H3) a control $v\in
L^{\infty }_{loc}(\overline{\omega }\times [0,+\infty ))$  exists
such that the solution $(u_1,u_2)$  satisfies
$$u_1(x,t)\geq 0, \quad u_2(x,t)\geq 0, \qquad \mbox{\rm a.e. } x\in \Omega , \ \mbox{\rm for any } t\geq 0$$
and
$$\displaystyle
\lim_{t\to \infty} \Vert u_1(t)\Vert_{L^{\infty }(\Omega)}=
\lim_{t\to \infty} \Vert u_2(t)\Vert_{L^{\infty }(\Omega)}=0 .$$

\bf Definition 3.2. \rm We say that our system  is \emph{locally
zero-stabilizable} if   there exists $r_0>0$ such that for any
$u_1^0$ and $u_2^0$ satisfying (H3) and $\|u_1^0\|_{L^{\infty
}(\Omega )}, \ \|u_2^0\|_{L^{\infty }(\Omega )}\leq r_0,$ there
exists $v\in L^{\infty }_{loc}(\overline{\omega }\times[0,+\infty
))$ such that the solution $(u_1,u_2)$  satisfies $u_1(x,t)\geq 0,
\ u_2(x,t)\geq 0$ a.e. $x\in \Omega , \ \mbox{\rm for any } t\geq
0$ and $\lim_{t\rightarrow +\infty}\| u_1(t)\| _{L^{\infty
}(\Omega )}= \lim_{t\rightarrow +\infty}\| u_2(t)\| _{L^{\infty
}(\Omega )}=0 .$

\bf Remark 2. \rm It is obvious that if a system is zero-stabilizable, then it
is also locally zero-stabilizable.

A stabilization result for our system, in the case of time
independent $g$, had been obtained in \cite{AK02}.
  In case of stabilizability a complicated stabilizing control had been
provided. A stronger result (which indicates also a simpler
stabilizing control) has been established   in \cite{AK09}  using
a different approach.  Later, in \cite{AKBelleni}  the authors
have further extended the main results to the case of a time
$T-$periodic function $g$ and provided a very simple stabilizing
feedback control.

In  \cite{AK09}, by  Krein-Rutman Theorem,  it has been shown that

$$\left\{ \begin{array}{ll}
- \displaystyle d_1\Delta \varphi +a_{11}\varphi
-{a_{21}\over{a_{22}}}\int_{\Omega }k(x,x') \varphi (x')dx'
=\lambda \varphi ,  \ &x\in \Omega \setminus \overline{\omega }\\
\varphi (x)=0,  &x\in \partial \omega  \\ \displaystyle {{\partial
\varphi }\over {\partial \nu }}(x)+\alpha \varphi (x)=0,&x\in
\partial \Omega ,
          \end{array}
  \right. $$
   admits a principal
(and real) eigenvalue $\lambda_{1 }(\omega),$ and a corresponding
strictly positive eigenvector $\varphi \in Int (\mathcal K)$ where

$${\mathcal K}=\{ \varphi \in L^{\infty }(\Omega); \ \varphi (x)\geq 0 \ \mbox{\rm a.e. in } \Omega\} . $$

The following theorem holds  \cite{AK09}:

\bf Theorem 3.3. \it If $\lambda _1(\omega )>0$, then for $\gamma
\geq 0$ large enough, the feedback control $v:=-\gamma u_1$
stabilizes our system to zero.

Conversely, if $h$ is differentiable at $0$ and $h'(0)=a_{21}$ and
if our system  is zero-stabilizable, then $\lambda_1(\omega )\geq
0$.
\rm
\vspace{3mm}

Moreover, the proof of the main result in \cite{AK09} shows that
for a given affordable sanitation effort $\gamma $, the epidemic
process can be diminished exponentially if $\lambda _{1,\gamma
}^{\omega }>0$ (at the rate of $\exp \{-\lambda _{1,\gamma }^{\omega }t\}$) , where $\lambda _{1,\gamma }^{\omega }$ is the
principal eigenvalue to the following problem:
\begin{equation}\left\{
\begin{array}{ll} - \displaystyle d_1\Delta \varphi +a_{11}\varphi
-{a_{21}\over{a_{22}}}\int_{\Omega }k(x,x') \varphi (x')dx'
+\gamma \chi _{\omega }\varphi =\lambda \varphi ,
 \ &x\in \Omega \\
\displaystyle {{\partial \varphi }\over {\partial \nu }}(x)+\alpha
\varphi (x)=0,&x\in \partial \Omega .
          \end{array}
  \right. \label{(1.2)}\end{equation}

A natural question related to the practical implementation of the
sanitation policy is the following: ``For a given sanitation
effort $\gamma >0$ in the region $\omega $, is the principal
eigenvalue $\lambda _{1,\gamma }^{\omega }$ positive (and
consequently can our  epidemic system  be stabilized to zero by
the feedback control $v:=-\gamma u_1$) ?''

So, the first problem to be treated  is the estimation of $\lambda
_{1,\gamma }^{\omega }$. Since  this eigenvalue problem is related
to a non-self adjoint
 operator, we cannot use a variational principle (as Rayleigh's for
selfadjoint operators); hence in \cite{AK12} the authors  have
proposed an alternative method  based  on the    following result:

\begin{equation}\lim_{t\rightarrow +\infty }\int_{\Omega }y^{\omega
}(x,t)dx=\zeta -\lambda _{1,\gamma }^{\omega } ,\label{(1.3)}
\end{equation}
where $y^{\omega }$ is the unique positive solution to

\begin{equation}\left\{
\begin{array}{ll} \displaystyle {{\partial y}\over {\partial t}}-
d_1\Delta y+a_{11}y+\gamma \chi_{\omega }y-{{a_{21}}\over
{a_{22}}}
\int_{\Omega }k(x,x') y(x',t)dx' & \\
~~~~~~~~~~~~~~~~~~\displaystyle -\zeta y+(\int_{\Omega }y(x,t)dx)y=0, \quad &x\in \Omega , \ t>0 \\
\displaystyle {{\partial y}\over {\partial \nu }}(x,t)+\alpha y(x,t)=0, \quad &x\in \partial \Omega , \ t>0 \\
y(x,0)=1, \quad &x\in \Omega ,
          \end{array}
  \right. \label{(1.4)}\end{equation}
and $\zeta >\lambda _{1,\gamma }^{\omega }$ is a constant.

\bf Remark 4. \rm Problem (\ref{(1.4)}) is a logistic model for the population
dynamics with diffusion and migration.   Since the solutions to
the logistic models rapidly stabilize, this means that
(\ref{(1.3)}) gives an efficient method to approximate $\lambda
_{1,\gamma }^{\omega }$. Namely, for $T>0$ large enough,
$$\zeta -\int_{\Omega }y^{\omega }(x,T)dx$$
gives a very good approximation of $\lambda _{1,\gamma }^{\omega
}$. The  above result leads  to a concrete numerical estimation of
$\lambda _{1,\gamma }^{\omega }$ by  analyzing the large-time
behavior  of the system for  different  values of $\zeta.$

 We  may also  remark that, if in   (\ref{(1.4)})
\begin{equation}y(x,0)=y_0, \quad x\in \Omega ,\label{(999)} \end{equation} with $y_0$ an
arbitrary positive constant, then
$$\lim_{t\rightarrow +\infty }\int_{\Omega }y_1^{\omega }(x,t)dx=\zeta -\lambda _{1,\gamma }^{\omega } ,$$
where $y_1^{\omega }$ is the solution to
(\ref{(1.4)})-(\ref{(999)}).

Assume now that for a given sanitation effort $\gamma $, the
principal eigenvalue to (\ref{(1.2)}) satisfies $\lambda
_{1,\gamma }^{\omega }>0$, and consequently $v:=-\gamma u_1$
stabilizes to zero the solution to (\ref{(1.1)}).

Let $\omega _0$ be a nonempty open subset of $\Omega $, with a
smooth boundary and such that $\omega _0\subset \subset \Omega $
and $\Omega \setminus \overline{\omega }_0$ is a domain. Consider
${\mathcal O}$ the set of all translations $\omega $ of $\omega
_0$, satisfying $\omega \subset \subset \Omega $. Since, after
all, our initial goal was to eradicate the epidemics, we are led
to the natural problem of ``Finding the translation $\omega ^*$ of
$\omega $ ($\omega \in {\mathcal O}$) which gives a small value
(possibly minimal) of
$$R^{\omega }=\int_{\omega }[u_1^{\omega }(x,T)+u_2^{\omega
}(x,T)]dx ,$$  at some given finite time $T>0.$    ''

Here $(u_1^{\omega}, u_2^{\omega})$ is the solution of (1.1)
corresponding to $v:=-\gamma u_1$, i.e. $(u_1^{\omega},
u_2^{\omega})$ is the solution to

\begin{equation}\left\{ \begin{array}{ll} \displaystyle {{\partial u_1}\over
{\partial t}}(x,t)= \displaystyle d_1\Delta
u_1(x,t)-a_{11}u_1(x,t)
+\int_{\Omega }&k(x,x') u_2(x',t)dx' \vspace{2mm} \\
\qquad \qquad \quad -\gamma \chi _{\omega }(x) u_1(x,t),\quad
&(x,t)\in
\Omega \times (0,+\infty ) \\
\displaystyle {{\partial u_1}\over {\partial \nu }}(x,t)+\alpha
u_1(x,t)=0,
\quad &(x,t)\in \partial \Omega \times (0,+\infty ) \\
\displaystyle {{\partial u_2}\over {\partial
t}}(x,t)=-a_{22}u_2(x,t) +g(u_1(x,t)), \ \
&(x,t)\in \Omega \times (0,+\infty ) \\
u_1(x,0)=u_1^0(x), &x\in \Omega \\
u_2(x,0)=u_2^0(x), &x\in \Omega .
          \end{array}
  \right. \label{(1.6)}
  \end{equation}

For this reason we are going to evaluate the derivative of
$R^{\omega }$ with respect to translations of $\omega $. This will
allow to derive a conceptual iterative algorithm to improve at
each step the position (by translation) of $\omega $ in order to
get a smaller value for $R^{\omega }$.

\subsection{ The derivative of $R^{\omega }$ with respect to
translations}

For any $\omega \in {\mathcal O}$ and $V\in \mbox{\bf R}^n$ we
define the derivative
$$dR^{\omega }(V)=\lim_{\varepsilon \rightarrow 0}{{R^{\varepsilon V+\omega }-R^{\omega }}\over {\varepsilon
}}.$$ For basic results and methods in the optimal shape design
theory we refer to \cite{Henrot}.

\bf Theorem 3.4. \label{3.1} \it For any $\omega \in {\mathcal O}$ and $V\in \mbox{\bf R}^n$ we
have that
$$dR^{\omega }(V)=\gamma \int_0^T\int_{\partial \omega }u_1^{\omega }(x,t)p_1^{\omega }(x,t)\nu (x)\cdot Vd\sigma \ dt, $$ where
$(p_1^{\omega },p_2^{\omega })$ is the solution to the adjoint
problem

\begin{equation}\left\{ \begin{array}{ll} \displaystyle {{\partial p_1}\over
{\partial t}}+d_1\Delta
p_1-a_{11}p_1-\gamma \chi_{\omega }p_1+g'(u_1^{\omega })p_2=0, & x\in \Omega , \ t>0 \\
\displaystyle {{\partial p_2}\over {\partial t}} + \int_{\Omega
}k(x',x)p_1(x',t)dx'-a_{22}p_2 =0, & x\in \Omega , \ t>0 \\
\displaystyle {{\partial p_1}\over {\partial \nu }}(x,t)+\alpha p(x,t)=0, \quad &x\in \partial \Omega , \ t>0 \\
p_1(x,T)=p_2(x,T)=1, \quad &x\in \Omega .
          \end{array}
  \right. \label{(3.1)}
  \end{equation}
Here $\nu (x)$ is the normal inward versor at $x\in \partial
\omega $ (inward with respect to $\omega $).
\rm

For the construction of the adjoint problems in optimal control
theory we refer to \cite{Lions_1988}.

 Based on Theorem \ref{3.1}, in  \cite{AK12}  the authors have proposed
 a conceptual iterative algorithm to improve the position (by
translation) of $\omega \in {\mathcal O}$ (in order to obtain a
smaller value for $R^{\omega }$.

\subsection{The periodic case}

As a purely technical simplification, we have  assumed that only
the incidence rate is periodic, and in particular that it can be
expressed as

$$(i.r.)(x,t)=h(t,u_1(x,t))=p(t)g(u_1(x,t)),$$
were $g$, the functional dependence of the incidence rate upon the
concentration of the pollutant, can be chosen as in the time
homogeneous case.

 In this case our goal is to study the controlled
system
\begin{equation}\left\{
\begin{array}{ll}
\displaystyle{\frac{\partial u_1 }{\partial t }}(x,t)=
\displaystyle d_1\Delta u_1(x,t)-a_{11}u_1(x,t) +\int_{\Omega }
&k(x,x') u_2(x',t)dx'\\ \quad \qquad \qquad
 + \chi_{\omega}(x) v(x,t), &(x,t)\in \Omega \times (0,+\infty ) \\
\displaystyle{\frac{\partial u_1 }{\partial \nu }}(x,t)+\alpha
u_1(x,t)=0,
\quad &(x,t)\in \partial \Omega \times (0,+\infty ) \\
\displaystyle{\frac{\partial u_2 }{\partial t
}}(x,t)=-a_{22}u_2(x,t) +h(t,u_1(x,t)), \ \
&(x,t)\in \Omega \times (0,+\infty ) \\
u_1(x,0)=u_1^0(x), \quad u_2(x,0)=u_2^0(x), &x\in \Omega ,
          \end{array}
  \right. \label{1.2}\end{equation}
with a control $v\in L^{\infty }_{loc}(\overline{\omega }\times
[0,+\infty ))$ (which implies that $supp (v(\cdot,t))\subset
{\overline \omega}$ for $t\geq 0$).

The explicit time dependence of the incidence rate is given via
the function $p(\cdot),$  which is assumed to be a strictly
positive, continuous and $T-$periodic function of time; i.e. for
any  $t \in \mathbb R,$

$$p(t)=p(t+T).$$

\bf Remark 4. \rm The  results  can be easily extended to the case
in which  also $a_{11}$, $a_{22}$ and $k$ are $T-$periodic
functions.

Consider  the following (linear) eigenvalue problem
\begin{equation}\left\{ \begin{array}{ll}
\displaystyle{\frac{\partial \varphi }{\partial t}}  -
\displaystyle d_1\Delta \varphi +a_{11}\varphi -\int_{\Omega
}k(x,x') \psi(x',t)dx' =\lambda \varphi ,
&x\in \Omega \setminus \overline{\omega }, \ t>0\\
\displaystyle{\frac{\partial \varphi }{\partial \nu}} (x,t)+\alpha
\varphi (x,t)=0,&x\in \partial \Omega , \ t>0\\
\varphi (x,t)=0, &x\in \partial \omega , \ t>0\\
\displaystyle{\frac{\partial \psi }{\partial t}} (x,t)+a_{22}\psi
(x,t)-a_{21}p(t)\varphi
(x,t)=0,\quad &x\in \Omega \setminus \overline{\omega }, \ t>0\\
\varphi (x,t)=\varphi (x,t+T), \quad \psi (x,t)=\psi (x,t+T),
\quad &x\in \Omega \setminus \overline{\omega }, \ t\geq 0.
          \end{array}
  \right. \label{1.3}\end{equation}
By similar procedures,  as  in the time homogeneous case, Problem
(\ref{1.3}) admits a principal (real) eigenvalue $\lambda _{1
}^T(\omega),$ and a corresponding strictly positive eigenvector
$\varphi^T \in Int (\mathcal K^T)$ where
$${\mathcal K^T}=\{ \varphi \in L^{\infty }(\Omega \times (0,T)); \ \varphi (x,t)\geq 0 \ \mbox{\rm a.e. in } \Omega \times (0,T)\} . $$

\bf Theorem 3.5. \label{theo1} \it If $\lambda _1^T(\omega )>0$, then
for $\gamma \geq 0$ large enough, the feedback control $v:=-\gamma
u_1$ stabilizes (\ref{1.2}) to zero.

Conversely, if $g$ is differentiable at $0$ and $g'(0)=a_{21},$
and if (\ref{1.2}) is zero-stabilizable, then $\lambda _1^T(\omega
)\geq 0$.\rm
\vspace{3mm}

\bf Theorem 3.6 \label{theo2} \it Assume that $g$ is differentiable at
$0$. Denote by $\tilde{\lambda }_1^T(\omega )$ the principal
eigenvalue of the problem
\begin{equation}
\left\{ \begin{array}{ll} \displaystyle{\frac{\partial \varphi
}{\partial t}}  - \displaystyle d_1\Delta \varphi +a_{11}\varphi
-\int_{\Omega }k(x,x') \psi(x',t)dx'  =\lambda \varphi ,
&x\in \Omega \setminus \overline{\omega }, \ t>0\\
\displaystyle{\frac{\partial \varphi }{\partial \nu}}
(x,t)+\alpha
\varphi (x,t)=0,&x\in \partial \Omega , \ t>0\\
\varphi (x,t)=0, &x\in \partial \omega , \ t>0\\
\displaystyle{\frac{\partial \psi }{\partial t}} (x,t)+a_{22}\psi
(x,t)-g'(0)p(t)\varphi
(x,t)=0,\quad &x\in \Omega \setminus \overline{\omega }, \ t>0\\
\varphi (x,t)=\varphi (x,t+T), \quad \psi (x,t)=\psi (x,t+T),
\quad &x\in \Omega \setminus \overline{\omega }, \ t\geq 0
          \end{array}
  \right.\label{local}
  \end{equation}
If $\tilde{\lambda }_1^T(\omega )>0$, then the system is locally
zero stabilizable, and for $\gamma \geq 0$ sufficiently large,
$v:=-\gamma u_1$ is a stabilizing feedback control.

Conversely, if the system is locally zero stabilizable, then
$\tilde{\lambda }_1^T(\omega )\geq 0$.
\rm

 \vspace{1mm}

\bf Remark 5. \rm Since $g'(0)\leq a_{21}$, it follows that $\lambda _1^T(\omega )
\leq \tilde{\lambda }_1^T(\omega ).$ We conclude now that

\begin{itemize}

\item[$1^0$] If $\lambda _1^T(\omega )>0$, the system is
zero-stabilizable;

\item[$2^0$] If $\tilde{\lambda }_1^T(\omega )>0$ and $\lambda
_1^T(\omega )\leq 0$, the system is locally zero-stabilizable;

\item[$3^0$] If $\tilde{\lambda }_1^T(\omega )<0$, the system is
not locally zero-stabilizable and consequently it is not zero
stabilizable.
\end{itemize}

\bf Remark 6. \rm \emph{Future directions.} Another interesting problem is that when
$\omega $ consists of a finite number of mutually disjoint
subdomains. The goal is to find the best position for each
subdomain. A similar approach can be used.

In a  recently submitted  paper \cite{AK15}, the problem of the
best choice of  the subregion $\omega$ has been faced for a
general harvesting problem in population dynamics as a shape
optimization problem; our  future aim is to apply  those results
to  our problem of eradication of   spatially  structured
epidemics.

\section*{Acknowledgments}
It is a pleasure to acknowledge  the contribution  by  Klaus Dietz
regarding the  bibliography  on the historical remarks reported in
Section 1.



\begin{thebibliography}{99}




\bibitem{Abbey52} H. Abbey, \emph{An examination of the Reed-Frost theory of epidemics}, Human Biol. \textbf{24}
(1952), 201--233.

\bibitem{AK02} S. Ani\c{t}a, V. Capasso, \emph{A stabilizability problem for a reaction-diffusion system modelling a class of
spatially structured epidemic systems}, Nonlin. Anal. Real World Appl. \textbf{3} (2002), 453--464.

\bibitem{AK09} S. Ani\c{t}a, V. Capasso, \emph{A stabilization strategy for a reaction-diffusion system modelling a class of
spatially structured epidemic systems (\emph{think globally, act locally})}, Nonlin. Anal. Real World Appl. \textbf{10} (2009), 2026--2035.

\bibitem{AKBelleni} S. Ani\c{t}a, V. Capasso, \emph{On the stabilization of reaction-diffusion systems modeling a class of
man-environment epidemics: a review}, Math. Meth. Appl. Sci. \textbf{33} (2010), 1235--1244.

\bibitem{AKB} S. Ani\c{t}a, V. Capasso, \emph{Stabilization for a reaction-diffusion system modelling a class of spatially
structured epidemic systems. The periodic case}, in ``Advances in Dynamics and Control: Theory, Methods, and Applications'' (eds. S.
Sivasundaram, et al.),  Cambridge Scientific Publishers Ltd., Cambridge, MA, 2011.

\bibitem{AK12} S. Ani\c{t}a, V.  Capasso, \emph{Stabilization of  a reaction-diffusion system modelling
a class of spatially structured epidemic systems via feedback control}, Nonlin. Anal. Real World Appl. \textbf{13} (2012), 725--735.

\bibitem{AK15} S.  Ani\c{t}a, V. Capasso, A.-M. Mo\c{s}neagu, \emph{Regional control in optimal harvesting of population dynamics},
Nonlin. Anal. \textbf{147} (2016), 191--212.

\bibitem{aris_75} R. Aris, The Mathematical Theory of Diffusion and Reaction in Permeable Catalysts. Vols. I and II, Oxford Univ.
Press, Oxford, 1975.

\bibitem{VK_Barbu_Arnautu} V. Arn\u{a}utu, V. Barbu, V.  Capasso, \emph{Controlling the spread of a class of epidemics},
Appl. Math. Optimiz. \textbf{20} (1989), 297--317.

\bibitem{AronMay} J.L. Aron, R.M. May, \emph{The population dynamics of malaria}, in ``Population Dynamics of Infectious Diseases. Theory and Applications'' (ed. R.M. Anderson), Chapman
\& Hall, London, 1982, 139--179.

\bibitem{Aronson_1977} D.G. Aronson, \emph{The asymptotic speed of propagation of a simple epidemic},
in ``Nonlinear Diffusion'' (eds. W. E. Fitzgibbon, A.F. Walker), Pitman, London, 1977, 1--23.

\bibitem{aronson_85} D.G. Aronson, \emph{The role of diffusion in mathematical population biology : Skellam revisited}, in ``Mathematics in Biology and Medicine''
(eds.  V. Capasso, E. Grosso, S.L. Paveri-Fontana), Lecture Notes in Biomathematics, vol. 57, Springer-Verlag,
Berlin, Heidelberg, 1985.

\bibitem{aronson_weinberger} D.G. Aronson, H.F. Weinberger, \emph{Nonlinear diffusions in
population genetics, combustion and nerve propagation}, in ``Partial Differential Equations and Related
Topics'' (ed. J. Goldstein), Lecture Notes in Mathematics, vol. 446, Springer-Verlag,
Berlin, Heidelberg, 1975.

\bibitem{babak_2007} P. Babak, \emph{Nonlocal initial problems for coupled reaction-diffusion systems and their applications},  Nonlin. Anal. Real World Appl. \textbf{8} (2007), 980--996.

\bibitem{Bac} N. Bacaer, C. Sokhna, \emph{A reaction-diffusion system modeling the spread of resistance to an
antimalarial drug}, Math. Biosci. Eng. \textbf{2} (2005), 227--238.

\bibitem{Bailey50} N.T.J. Bailey,  \emph{A simple stochastic epidemic}, Biometrika \textbf{37} (1950), 193--202.

\bibitem{Bartlett49} M.S. Bartlett, \emph{Some  evolutionary stochastic processes}, J. Roy. Stat. Soc. Ser. B \textbf{11} (1949), 211--229.

\bibitem{Bernoulli1760} D. Bernoulli, \emph{R\'{e}flexions sur les avantages de l'inoculation}, Mercure de France \textbf{June} (1760), 173-190. Reprinted in \emph{ Die
Werke von Daniel Bernoulli}, Bd. 2 Analysis und Wahrscheinlichkeitsrechnung (eds. L.P. Bouckaert, B.L. van der Waerden), Birkh\"{a}user, Basel, 1982, p. 268.

\bibitem{Brownlee11} J. Brownlee,  \emph{The mathematical theory of random migration and epidemic distribution}, Proc. Roy. Soc. Edinburgh \textbf{31}
(1911), 262--288.

\bibitem{Busenberg_etal_1988} S. Busenberg, K. L. Cooke, M. Iannelli, \emph{Endemic thresholds and stability in a class of age-structured epidemics},
SIAM J. Appl. Math. \textbf{48} (1988), 1379--1395.

\bibitem{VK78} V. Capasso, \emph{Global solution for a diffusive nonlinear deterministic epidemic model}, SIAM J. Appl. Math. \textbf{35} (1978), 274--284.

\bibitem{VK84} V. Capasso, \emph{Asymptotic stability for an integro-differential reaction-diffusion system}, J. Math. Anal. Appl. \textbf{103} (1984), 575--588.

\bibitem{VK09} V. Capasso, Mathematical Structures of Epidemic Systems (2nd revised printing), Lecture Notes Biomath., vol. ~97, Springer-Verlag, Heidelberg, 2009.

\bibitem{VK83} V. Capasso, L. Maddalena, \emph{Periodic solutions for a reaction-diffusion system modelling the spread of a
class of epidemics}, SIAM J. Appl. Math. \textbf{43} (1983), 417--427.

\bibitem{7} V. Capasso, L. Maddalena, \emph{On a degenerate nonlinear diffusion problem with boundary feedback}, Appl. Anal. \textbf{24} (1987), 256--298.


\bibitem{VK_Kunisch} V. Capasso, K. Kunisch, \emph{A reaction-diffusion system arising in modelling man-environment diseases}, Quarterly Appl. Math. \textbf{46}
(1988), 431--450.

\bibitem{VK81} V. Capasso, L. Maddalena, \emph{Convergence to equilibrium states for a reaction-diffusion system
modelling the spatial spread of a class of bacterial and viral diseases}, J. Math. Biol. \textbf{13} (1981), 173--184.

\bibitem{VK82} V. Capasso, L. Maddalena, \emph{Saddle point behaviour for a reaction-diffusion system: Application to a class of epidemic models}, Math. Comput. Simulation \textbf{24} (1982), 540--547.

\bibitem{VK83} V. Capasso, L.  Maddalena, \emph{Periodic solutions for a reaction-diffusion system modelling the spread of a class of epidemics}, SIAM J. Appl. Math. \textbf{43} (1983),  417--427.

\bibitem{18.} V. Capasso,  S.L. Paveri-Fontana, \emph{A mathematical model for the 1973 cholera epidemic in the European Mediterranean
region}, Revue d'Epidemiologie et de la Sant\'{e} Publique \textbf{27} (1979), 121--132. Errata corrig\'{e} \textbf{28} (1980), 390.

\bibitem{12.} V. Capasso, G. Serio, \emph{A generalization of the Kermack-Mckendrick deterministic epidemic model}, Math. Biosci. \textbf{42} (1978), 43--61.

\bibitem{VK_Thieme} V.  Capasso, H.R. Thieme, \emph{A threshold theorem for a reaction-diffusion epidemic system}, in ``Differential Equations and
Applications'' (ed. R. Aftabizadeh), Ohio Univ. Press, Athens, OH, 1989, pp. 128--133.

\bibitem{VK97} V. Capasso, R.E.  Wilson, \emph{Analysis of a reaction-diffusion system modelling man-environment-man epidemics}, SIAM J. Appl. Math. \textbf{57} (1997), 327--346.

\bibitem{Castillo_Chavez_1_1989} C. Castillo-Chavez, K.L. Cooke, W. Huang, S.A. Levin, \emph{On the  role of long incubation periods in the dynamics of acquired immunodeficiency syndrome
(AIDS). Part 1: Single  population models}, J. Math. Biol. \textbf{27} (1989), 373--398.

\bibitem{Castillo_Chavez_2_1989} C. Castillo-Chavez, K.L. Cooke, W. Huang, S.A. Levin, \emph{On the  role of long incubation periods in the dynamics of acquired immunodeficiency syndrome
(AIDS). Part 2: Multiple group models}, in ``Mathematical and Statistical Approaches to  AIDS Epidemiology''  (ed. C. Castillo-Chavez), Lecture Notes in Biomathematics, Vol. 83, Springer-Verlag, Heidelberg, 1989, 200--217.

\bibitem{britton} F. Chamchod, N. F. Britton, \emph{Analysis of a vector-bias model on malaria transmission},
Bull. Math. Biol. \textbf{73} (2011), 639--657.

\bibitem{Chitnis_2006} N. Chitnis, J.M. Cushing, J. M. Hyman, \emph{Bifurcation analysis of a mathematical model for
malaria transmission}, SIAM J. Appl. Math. \textbf{67} (2006), 24--45.

\bibitem{Chitnis_2008} N. Chitnis, T.A. Smith, R.W. Steketee, \emph{A mathematical model for the dynamics of malaria in
mosquitoes feeding on a heterogeneous host population}, J. Biol. Dyn. \textbf{2} (2008), 259--285.

\bibitem{cliff_ord} A.D. Cliff, P. Haggett, J.K. Ord, G.R. Versey, Spatial Diffusion. An Historical Geography of Epidemics in an Island
Community, Cambridge University Press, Cambridge, 1981.

\bibitem{Cod} C.T. Code\c{c}o, \emph{Endemic and epidemic dynamics of cholera: the role of the aquatic reservoir}, BMC Infectious Diseases \textbf{1} (2001), 1--14.

\bibitem{diekmann_78} O. Diekmann, \emph{Thresholds and travelling waves for the geographical spread of infection}, J. Math. Biol. \textbf{6} (1978),
109--130.

\bibitem{Dietz_1} K. Dietz, \emph{Mathematical models for transmission and control of malaria},  in ``Principles and Practice of Malariology'' (eds. W. Wernsdorfer, Y. McGregor), Churchill
Livingstone, Edinburgh, 1988, 1091--1133.

\bibitem{Dietz_1997} K. Dietz, \emph{Introduction to McKendrick (1926) applications of mathematics to medical problems},
in ``Breakthroughs in Statistics'', Volume III,  eds. S. Kotz, N.L. Johnson), Springer-Verlag, Heidelberg, 1997, 17--26.

\bibitem{Dietz_2009} K. Dietz, \emph{Mathematization in sciences epidemics: the fitting of the first dynamic models to data}, J. Contemp. Math. Anal. \textbf{44}
(2009), 36--46.

\bibitem{Dietz} K. Dietz, T. Molineaux, A. Thomas, \emph{A malaria model tested in the African savannah},
Bull. World Health Org. \textbf{50} (1974), 347--357.

\bibitem{DOnofrio} A. d'Onofrio, P. Manfredi, E. Salinelli, \emph{Vaccination behaviour, information, and the dynamics of SIR vaccine preventable diseases},
 Theor. Popul. Biol. \textbf{71} (2007), 301--317.

\bibitem{Enko1989} En'ko, \emph{On the course of epidemics of some infectious diseases}, Translation from
Russian by K. Dietz, Int. J. Epidem. \textbf{18} (1989), 749--755.

\bibitem{Farr1840} W. Farr,  \emph{Progress of epidemics}, Second Report of the Registrar General (1840), 91--98.

\bibitem{fife} P.C. Fife, Mathematical Aspects of Reacting and Diffusing Systems,  Lecture Notes in Biomathematics, vol. 28,
Springer-Verlag, Berlin, Heidelberg, 1979.

\bibitem{Frost_1976} W.H. Frost, \emph{Some conceptions of epidemics in general}, Am. J. Epidem. \textbf{103} (1976), 141--151.

\bibitem{Hamer06} W. H. Hamer, \emph{Epidemic disease in England}, Lancet \textbf{1} (1906), 733--739.

\bibitem{Henrot} A. Henrot, M. Pierre, Variation et Optimisation de Formes. Une Analyse G\'{e}om\'{e}trique,
Springer, Berlin, 2005.

\bibitem{Hethcote_vandendriessche_91} H.W. Hethcote, P. van den Driessche, \emph{Some epidemiological models with  nonlinear incidence},
J. Math. Biol. \textbf{29} (1991), 271--287.

\bibitem{Hoppensteadt_75} F. Hoppensteadt, Mathematical Theories  of Populations: Demographics, Genetics and Epidemics, SIAM, Philadelphia,
1975.

\bibitem{Kendall_57} D.G. Kendall, \emph{Mathematical models of the spread of infection}, in ``Mathematics and Computer Science in
Biology and Medicine'', H.M.S.O., London, 1965, 213--225.

\bibitem{Kermack_McKendrick27} W.O. Kermack, and A. G. McKendrick, \emph{A contribution to the mathematical theory of epidemics},
Proc. Roy. Soc. London, Ser. A \textbf{115} (1927), 700--721.

\bibitem{Killeen} G.F. Killeen, F.F. McKenzie, B.D. Foy, C. Schieffelin, P.F. Billingsley, J.C. Beier,
\emph{A simplified model for predicting malaria entomological inoculation rates based on entomologic
and parasitologic parameters relevant to control}, Am. J. Trop. Med. Hyg. \textbf{62} (2000), 535--544.

\bibitem{levin_78} S.A. Levin, \emph{Population models and community structure in heterogeneous environments}, in ``Mathematical
Association of America Study in Mathematical Biology. Vol. II: Populations and Communities'' (ed. S. Levin), Math. Assoc. Amer., Washington, 1978, 439--476.

\bibitem{Lions_1988} J. L. Lions, Controlabilit\'{e} exacte, stabilisation et perturbation de syst\`{e}mes distribu\'{e}s, RMA ~8, Masson, Paris, 1988.

\bibitem{Liu_Hethcote_Levin_1987} W.M. Liu, H.M. Hethcote, S.A. Levin, \emph{Dynamical behavior of epidemiological models with nonlinear incidence rates},
J. Math. Biol. \textbf{25} (1987), 359--380.

\bibitem{Lotka23} A. J. Lotka, \emph{Martini's  equations for the  epidemiology of immunizing diseases}, Nature \textbf{111} (1923), 633--634.

\bibitem{Macdonald50} G. Macdonald, \emph{The  analysis of malaria parasite rates in infants}, Tropical Disease Bull. \textbf{47}
(1950), 915--938.

\bibitem{mcdonald} G. Macdonald, The Epidemiology and Control of Malaria, Oxford University Press, London, 1957.

\bibitem{Martini_1921} E. Martini, Berechnungen und Beobachtungen zur Epidemiologie und Bek\"{a}mpfung der
Malaria, Gente, Hamburg, 1921.

\bibitem{molin} L. Molineaux, G. Gramiccia, The Garki Project. Research on the Epidemiology and Control of Malaria in the Sudan
Savanna of West Africa, World Health Organization, Geneva,1980.

\bibitem{murray_89} J.D. Murray, Mathematical Biology, Springer-Verlag, Berlin, Heidelberg, 1989.

\bibitem{nasell_hirsch} I. N\"{a}sell, W.M. Hirsch, \emph{The transmission dynamics of schistosomiasis}, Comm. Pure Appl. Math. \textbf{26} (1973), 395--453.

\bibitem{okubo} A. Okubo, Diffusion and Ecological Problems: Mathematical Models, Springer-Verlag,
Berlin, Heidelberg, 1980.

\bibitem{Puma_1939} M. Puma, Elementi per una Teoria Matematica del Contagio,
Editoriale  Aeronautica, Rome, 1939.

\bibitem{ross} R. Ross, The Prevention of Malaria, 2$^{nd}$ edition,  Murray, London, 1911.

\bibitem{ruan} S. Ruan, D. Xiaob, J. C. Beierc, \emph{On the delayed Ross-Macdonald model for malaria transmission},
 Bull. Math. Biol. \textbf{70} (2008), 1098--1114.

\bibitem{Serfling_1952} R.E. Serfling, \emph{Historical review of epidemic theory}, Human Biol. \textbf{24} (1952), 145--165.

\bibitem{Severo_1969} N.C. Severo, \emph{Generalizations  of  some stochastic epidemic models}, Math. Biosci. \textbf{31} (1969), 395--402.

\bibitem{skellam} J.G. Skellam, \emph{Random dispersal in theoretical populations}, Biometrika \textbf{38} (1951), 196--218.

\bibitem{Smith} T. Smith, G. Killeen, N. Maire, A. Ross, L. Molineaux, F. Tediosi, G. Hutton, J. Utzinger, K. Dietz, M. Tanner,
    \newblock \emph{Mathematical modeling of the impact of malaria vaccines on the clinical epidemiology and natural history of \textit{Plasmodium
Falciparum} malaria: overview}, Am. J. Trop. Med. Hyg. \textbf{75} (suppl. 2) (2006), 1--10.

\bibitem{smoller} J. Smoller, Shock Waves and Reaction-Diffusion Equations, Springer-Verlag, Berlin, Heidelberg, 1983.

\bibitem{Sochanta} T. Sochantha, S. Hewitt, C. Nguon, L. Okell, N. Alexander, S. Yeung, H. Vannara, M. Rowland, D. Socheat,
\emph{Insecticide-treated bednets for the prevention of Plasmodium falciparum malaria in Cambodia: a cluster-randomized trial},
Trop. Med. Int. Health \textbf{11} (2006), 1166--1177.

\bibitem{Soper29} H.E. Soper, \emph{The  interpretation of periodicity in disease prevalence}, J. Roy. Stat. Soc. \textbf{92}
(1929), 34--73.

\bibitem{thieme_77} H.R. Thieme, \emph{A model for the spatial spread of an epidemic}, J. Math. Biol. \textbf{4} (1977), 337--351.

\bibitem{WW_PNAS_1945} E.B. Wilson, J. Worcester, \emph{The  law of mass action in epidemiology}, Proc. Nat. Acad. Sci. \textbf{31} (1945), 24--34.

\bibitem{WHO_UNICEF} The Africa Malaria Report, WHO-UNICEF, 2003.

\end{thebibliography}
\end{document}